# On Flexible Prismatic Polyhedra

Gerald D. Nelson
nelso229@comcast.net

**Abstract**

We demonstrate the existence of four types of flexible prismatic polyhedra that can be derived or inferred from a consideration of Bricard octahedra and generalizations of Bricard octahedra. These flexible polyhedra are of genus 0 and 1, have dihedral angles all of which are non-constant under flexion, have faces that exhibit self-intersections, have an indefinite number of faces and have vertex motion that can treated as a function of a single variable.

## 1. Introduction

Prismatic polyhedra consisting of trapezoidal faces (quadrilaterals having two parallel edges) that form annular tubelike structures and that are flexible have been known for over a century, many of which are based upon the three types of flexible octahedra discovered in 1897 by Bricard [1]. An early paper [2] by Bennett published in 1911 describes a prismatic octahedron that is derived from a Bricard octahedron of the first type by letting opposite vertexes of the octahedron go to infinity. This result was expanded considerably in 1947 by Goldberg in a paper [3] in which necessary and sufficient conditions for flexion were developed for prismatic polyhedra having a quadrangle cross section and some of the resulting prismatic constructions were shown to be derivable from Bricard octahedra. A large flexible polyhedra described [4] in 1995 by Alexandrov for the purpose of providing a possible counter-example to the Bellows Conjecture is prismatic, although not explicitly based upon Bricard octahedra. Another that is based upon Bricard octahedra of the third type is mentioned briefly in a more recent paper [5] published by Stachel in 2004.

In this paper we describe four types of flexible prismatic polyhedra that can be constructed from a consideration of Bricard octahedra and of generalizations of Bricard octahedra as applied to flexible suspensions (polyhedra having the combinatorial structure of dipyramids).

We consider a flexible prismatic polyhedron to be a surface in $\mathbb{R}^3$ that is comprised of two or more segments each of which consists of simply connected trapezoidal faces joined together at their parallel edges to form an annular ring that is polygonal in cross section. These segments are joined at their ends at junctures whose geometric characteristics permit change in the dihedral angles at the junctures and along the parallel segment edges without any change in the trapezoidal faces. Generally some faces exhibit self-intersections, either at the juncture between annular segments or within the faces of the segments or between annular segments that are disjoint from one another. A genus 0 prismatic polyhedron has two or more segments the first and last segments of which are of infinite extent. A genus 1 prismatic polyhedron is a closed polyhedron that consists of four or more annular segments all of which are of finite extent and in which all segment ends are connected at junctures.

The approach that will be pursued in the following sections is to take a flexible suspension, as illustrated in Fig. 1, that has two apical vertexes **u**, **w** and N (an even integer ≥4) non-apical vertexes **$v_k$** for k=1..N and let the apical vertexes go to infinity as described in [2, Sec. 19] while keeping geometric relationships fixed at the non-apical vertexes; effectively letting the face angles $α_k$ and $A_k$ go to zero.

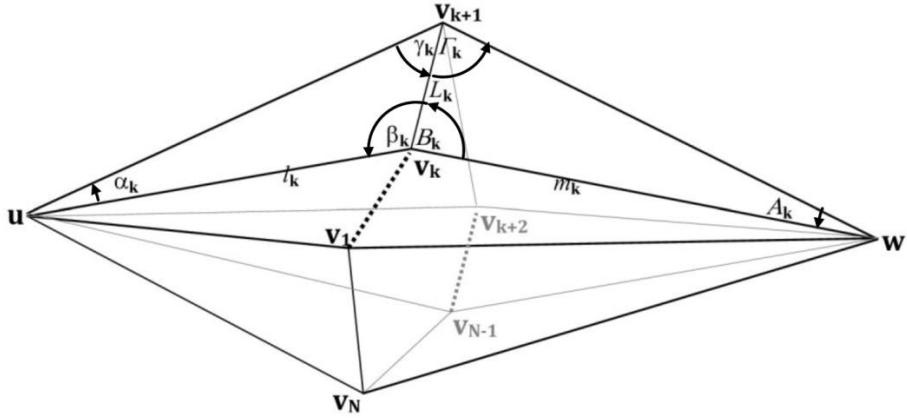

**Figure 1.** Polyhedral Suspension with apical vertexes of Index N

This results in the creation of an unbounded flexible suspension as illustrated in Fig. 2. In this manner we define the geometric characteristics, edge lengths and face angles, of the juncture between two annular segments so that the flexibility of the finite suspension (Fig. 1) is preserved.

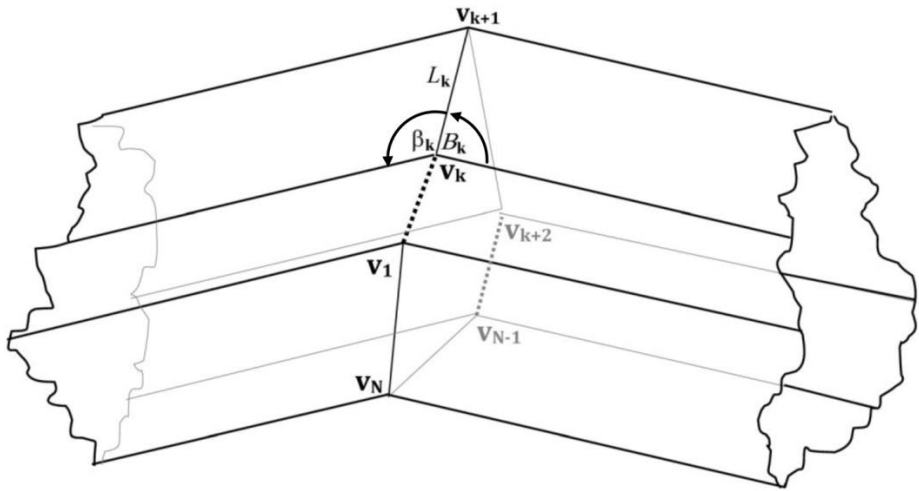

**Figure 2.** Unbounded Suspension with N vertexes

Additionally it is easy to see that the flexibility of the juncture that is characterized by the parameter set {$β_k$, $B_k$, $L_k$, k=1..N} as shown in Fig. 2 carries over into three alternate parameter sets that are created from reflecting the various trapezoidal segments at the juncture. These sets are the following: {$β_k$, $π-B_k$, $L_k$, k=1..N}, {$π-β_k$, $π-B_k$, $L_k$, k=1..N} and {$π-β_k$, $B_k$, $L_k$, k=1..N} as illustrated respectively by the alternate junctures shown in Figs. 3A, 3B and 3C.

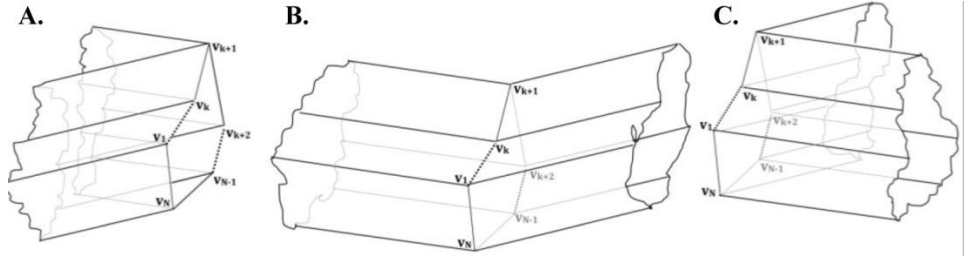

**Figure 3.** Alternate Suspension Junctures.
**A.** $\{\beta_k, \pi-B_k, L_k, k=1..N\}$  **B.** $\{\pi-\beta_k, \pi-B_k, L_k, k=1..N\}$  **C.** $\{\pi-\beta_k, B_k, L_k, k=1..N\}$

We proceed as follows: Sec.2 describes the characteristics of unbounded suspensions that are of interest. Sec. 3 provides arguments that these unbounded suspensions are flexible and a description of their construction. Sec. 4 describes how unbounded suspensions can be used to construct flexible prismatic polyhedra. Sec. 5 has a proof that the constructed prismatic polyhedra are flexible.

## 2. Unbounded Suspensions

Here we describe specific unbounded suspensions that are formed from flexible suspensions in which apical vertexes are allowed to go to infinity while geometric relationships at non-apical vertexes remain fixed. We limit our consideration to four types of flexible suspensions that have been previously defined [6]. These four types are designated by I-OEE, II-AEE, II-OEE and III-OAE. (A fifth type, III-OAS, is not of interest as the respective apical vertex angles exhibit sums which are $=\pi$ or $=2\pi$ and cannot all tend to zero.) In these names Roman numerals refer to the type of the related Bricard octahedron while the other letters are acronyms for geometric characteristics that can be associated with the type: OEE – Opposite Edges Equal, AEE – Adjacent Edges Equal, OAE – Opposite Angles Equal and OAS – Opposite Angles Supplementary. The four types of flexible suspensions and the results of their infinite elongation are described in the following using the notation illustrated in Figs. 1 and 2.

**I-OEE Suspensions** are generalizations of Bricard octahedra of the first type and are characterized by a total of 3M (where N=2M) independent edge length parameters out of a possible 3N parameters. These edge lengths have the following relationships:

$$\begin{aligned} m_k &= l_{k+M}, \\ m_{k+M} &= l_k \text{ and} \\ L_{k+M} &= L_k \text{ for } k = 1..M. \end{aligned} \quad (1)$$

As a consequence of equations (1) face angles at the non-apical vertexes of a suspension have the relationships:

$$\begin{aligned} B_k &= \beta_{k+M}, \\ \Gamma_k &= \gamma_{k+M}, \\ B_{k+M} &= \beta_k \text{ and} \\ \Gamma_{k+M} &= \gamma_k \text{ for } k = 1..M. \end{aligned} \quad (2)$$

By letting apical angles of the I-OEE suspension go to zero and retaining the angular and length relationships at the non-apical vertexes, equations (1) and (2) reduce to the following:

$$\begin{aligned} L_{k+M} &= L_k, \\ B_k &= \beta_{k+M} \text{ and} \\ B_{k+M} &= \beta_k \text{ for } k = 1..M. \end{aligned} \quad (3)$$

Also as apical angles go to zero segment edges become parallel and $\gamma_k = \pi - \beta_k$ and $\Gamma_k = \pi - B_k$ for all k. While equations (3) suggest that there are 3M independent parameters there is a "continuity" constraint that reduces this number by one. This constraint is formulated so that the sum of all offsets encountered by traversing around all segments from the first vertex and returning to the first vertex is zero. That is to say:

and
$$\sum_{k=1..N} L_k \cos \beta_k = 0$$
$$\sum_{k=1..N} L_k \cos B_k = 0. \quad (4)$$

By substituting from equations (3) it is easy to see that the two equations (4) reduce to the single equation:

$$\sum_{k=1..M} L_k (\cos \beta_k + \cos \beta_{k+M}) = 0. \quad (5)$$

Thus, one of the following 3M parameters $\{L_k, k = 1..M; \beta_k, k = 1..N\}$ can be determined from equation (5) while the remaining 3M-1 parameters are independent.

**II-AEE Suspensions** result from the first of two possible generalizations of Bricard octahedra of the second type and are characterized by 3M independent edge lengths that have the following relationships:

$$\begin{aligned} m_1 &= l_1 \text{ and} \\ m_k &= l_{N-k+2} \text{ for } k = 2..N \text{ and} \\ L_{k+M} &= L_{M-k+1} \text{ for } k = 1..M. \end{aligned} \quad (6)$$

From equations (6) it is seen that the face angles at the non-apical vertexes are related by:

$$\begin{aligned} \Gamma_k &= \beta_{N-k+1} \text{ and} \\ B_k &= \gamma_{N-k+1} \text{ for } k = 1..N. \end{aligned} \quad (7)$$

These angular relationships, equations (7), and the length relationships, equations (6), at the non-apical vertexes of the suspension are retained as apical angles go to zero, leaving:

$$\begin{aligned} L_{N-k+1} &= L_k \text{ for } k = 1..M \text{ and} \\ B_k &= \beta_{N-k+1} \text{ for } k = 1..N. \end{aligned} \quad (8)$$

The continuity equations (4) when combined with these II-AEE relationships reduce to one constraint:

$$\sum_{k=1..M} L_k (\cos \beta_k + \cos \beta_{N-k+1}) = 0. \tag{9}$$

Thus, one parameter of the 3M parameters $\{L_k, k = 1..M; \beta_k \text{ for } k = 1..N\}$ can be determined from equation (9) and the remaining 3M-1 are independent.

**II-OEE Suspensions** result from a second generalization of Bricard octahedra of the second type and are characterized by 3M independent edge lengths that have the following relationships:

$$\begin{aligned} l_{k+M} &= l_k, \\ m_{k+M} &= m_k \text{ and} \\ L_{k+M} &= L_k \text{ for } k = 1..M. \end{aligned} \tag{10}$$

Thus face angles at the non-apical vertexes that reflect equations (10) have the following relationships:

$$\begin{aligned} \beta_{k+M} &= \beta_k, \\ \gamma_{k+M} &= \gamma_k, \\ B_{k+M} &= B_k \text{ and} \\ \Gamma_{k+M} &= \Gamma_k \text{ for } k = 1..M. \end{aligned}$$

Retaining these angular relationships and the length relationships in equations (10) at the non-apical vertexes of the suspension and letting apical angles go to zero we are left with:

$$\begin{aligned} L_{k+M} &= L_k, \\ \beta_{k+M} &= \beta_k \text{ and} \\ B_{k+M} &= B_k \text{ for } k = 1..M. \end{aligned} \tag{11}$$

For II-OEE suspensions the continuity equations (4) are independent of one another and leave two equations when combined with these II-OEE relationships:

$$\sum_{k=1..M} L_k \cos \beta_k = 0$$

and

$$\sum_{k=1..M} L_k \cos B_k.$$

Thus the parameter set $\{L_k, \beta_k \text{ and } B_k \text{ for } k = 1..M\}$ has 3M-2 independent parameters.

**III-OAE Suspensions** result from generalizations of Bricard octahedra of the third type that are created from the parameter set $\{l_1, m_1, l_2, m_2, L_1, L_3, \ldots, L_{N-3}\}$ and that have well defined angular characteristics at the non-apical vertexes. There are N-2 vertexes having opposite face angles that are equal (OAE) while there are two vertexes in which the opposite face angles are supplementary (OAS). The OAS vertexes are, by convention, the vertex $\mathbf{v_1}$ and a second vertex $\mathbf{v_L}$ for an index L in the interval [3,N-2]. (In Bricard octahedra (N=4) the apical vertexes are also OAE.) Additionally, there are two flat foldings in which all vertexes are co-planar that are defined respectively by 1) $\delta_1 = \delta_L = 0$ and $\delta_k = \pi$ and by 2) $\delta_1 = \delta_L = \pi$ and $\delta_k = 0$ for k=2..L-1 and k=L+1..N where $\delta_k$ is the dihedral angle at the edge $\mathbf{uv_k}$.

As a consequence of the OAE and OAS relationships that are retained when apical angles go to zero the number of independent angles at the juncture at the non-apical vertexes is reduced to two, $\beta_1=\beta$ and $B_1=B$ say, and the remainder are defined as follows:

When L is even ($L=2k_x$):

$$\beta_{2k-1} = \beta \text{ and}$$
$$B_{2k-1} = B \text{ for } k = 1..k_x \text{ and}$$
$$\beta_{2k-1} = \pi - \beta \text{ and}$$
$$B_{2k-1} = \pi - B \text{ for } k = (k_x + 1)..M,$$

$$\beta_{2k} = \pi - B \text{ and}$$
$$B_{2k} = \pi - \beta \text{ for } k = 1..(k_x - 1) \text{ and}$$
$$\beta_{2k} = B \text{ and}$$
$$B_{2k} = \beta \text{ for } k = k_x..M,$$

When L is odd ($L=2k_x+1$): (12)

$$\beta_{2k-1} = \beta \text{ and}$$
$$B_{2k-1} = B \text{ for } k = 1..k_x \text{ and}$$
$$\beta_{2k-1} = \pi - \beta \text{ and}$$
$$B_{2k-1} = \pi - B \text{ for } k = (k_x + 1)..M,$$

$$\beta_{2k} = \pi - B \text{ and}$$
$$B_{2k} = \pi - \beta \text{ for } k = 1..k_x \text{ and}$$
$$\beta_{2k} = B \text{ and}$$
$$B_{2k} = \beta \text{ for } k = (k_x + 1)..M.$$

As a consequence of the flat foldings additional constraints are introduced. Effectively the cumulative width in a flat folding must be zero.

For the folding 1) $\delta_1=\delta_L=0$:

$$\sum_{k=1..(L-1)} L_k \sin \beta_k - \sum_{k=L..N} L_k \sin \beta_k = 0$$

and

$$\sum_{k=1..(L-1)} L_k \sin B_k - \sum_{k=L..N} L_k \sin B_k = 0$$

while for the folding 2) $\delta_1=\delta_L=\pi$:

$$\sum_{k=1..(L-1)} -1^{k+1} L_k \sin \beta_k - \sum_{k=L..N} -1^{k+1} L_k \sin \beta_k = 0$$

and

$$\sum_{k=1..(L-1)} -1^{k+1} L_k \sin B_k - \sum_{k=L..N} -1^{k+1} L_k \sin B_k = 0.$$

Using the relationships expressed in equations (12) the preceding equations reduce to the following two cases. When L ($L=2k_x$) is even:

$$\sum_{k=1..(k_x-1)} L_{2k} = \sum_{k=k_x..M} L_{2k},$$

and

$$\sum_{k=1..k_x} L_{2k-1} = \sum_{k=(k_x+1)..M} L_{2k-1},$$

when L ($L=2k_x+1$) is odd: (13)

$$\sum_{k=1..k_x} L_{2k} = \sum_{k=(k_x+1)..M} L_{2k},$$

and

$$\sum_{k=1..k_x} L_{2k-1} = \sum_{k=(k_x+1)..M} L_{2k-1},$$

It is easy to verify that the continuity equations (4) reduce to equations that are satisfied by equations (13) and do not impose further constraints upon the parameters. For example, when L is even equations (4) reduce to:

$$\cos\beta \left( \sum_{k=1..k_x} L_{2k-1} - \sum_{k=(k_x+1)..M} L_{2k-1} \right) = \cos B \left( \sum_{k=1..(k_x-1)} L_{2k} - \sum_{k=k_x..M} L_{2k} \right)$$

For all cases it is apparent that there are N independent parameters consisting of the angles β and $B$ and N-2 of the edge lengths $L_1$, $L_2$, $L_3$…. $L_N$ with the remaining two lengths determined by equations (13). A few constraint equations for edge lengths at the juncture are shown in Table I.

| N | L | Constraint Equations |
|---|---|---|
| 4 | 3 | $L_3=L_1$; $L_4=L_2$. |
| 6 | 3 | $L_1=L_3+L_5$; $L_2=L_4+L_6$. |
|   | 4 | $L_1+L_3=L_5$; $L_2=L_4+L_6$. |
| 8 | 3 | $L_1=L_3+L_5+L_7$; $L_2=L_4+L_6+L_8$. |
|   | 4 | $L_1+L_3=L_5+L_7$; $L_2=L_4+L_6+L_8$. |
|   | 5 | $L_1+L_3=L_5+L_7$; $L_2=L_4+L_6+L_8$. |
|   | 6 | $L_1+L_3+L_5=L_7$; $L_2+L_4+L_6=L_8$. |

Table I. Specific Length Constraints

## 3. Unbounded Flexible Suspensions

In this section we prove that unbounded prismatic suspensions, as shown in Fig. 2, having parameters that satisfy the conditions given by equations (3), (8), (11) or (12) and (13) are flexible and co-incidentally show how they can be constructed. This involves the specification of vectors that define the orientation of the parallel edges of the faces and of a recursive definition of coordinates that defines the positions of vertexes relative to one another. For each suspension type defined in Sec. 2, specific variations are introduced into the general coordinates as required to prove flexibility.

The orientation of segments relative to one another can be achieved by aligning parallel edges with unit vectors **u**(-sin θ, -cos θ, 0) and **w**(sin θ, -cos θ, 0) as shown in Fig. 4.

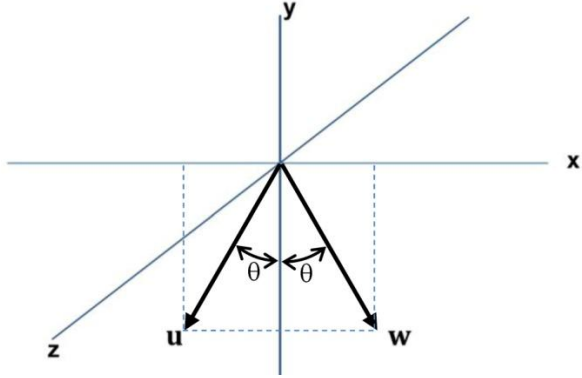

**Figure 4.** Geometry of orientation vectors.

Vertexes at the juncture are related by the equation:

$$\mathbf{v_{k+1}} = \mathbf{v_k} + \mathbf{\Delta v_k} \text{ for k} = 1..(N-1) \tag{14}$$

where the components of $\mathbf{\Delta v_k}$ are defined by:

$$\Delta x_k = \frac{L_k(\cos B_k - \cos \beta_k)}{2\cos \theta},$$

$$\Delta y_k = \frac{-L_k(\cos B_k + \cos \beta_k)}{2\cos \theta} \text{ and} \tag{15}$$

$$\Delta z_k = \pm\sqrt{L_k^2 - \Delta x_k^2 - \Delta y_k^2}.$$

The choice of the sign used for $\Delta z_k$ and the definition of $\mathbf{v_1}$ is type dependent. With the definitions as given by equations (15) we have by direct evaluation

$$|\mathbf{v_{k+1}} - \mathbf{v_k}| = L_k,$$
$$(\mathbf{v_{k+1}} - \mathbf{v_k})\bullet\mathbf{u} = L_k \cos \beta_k \text{ and}$$
$$(\mathbf{v_{k+1}} - \mathbf{v_k})\bullet\mathbf{w} = L_k \cos B_k \text{ for k} = 1..(N-1)$$

where ● is the vector inner product operator.

Since these results are only dependent upon edge lenghts and face angles it is clear that the angle θ can be treated as the variable of flexion and that equations (15) yield a flexible juncture provided that for the case when k=N we have:

$$\mathbf{v_1} = \mathbf{v_N} + \mathbf{\Delta v_N},$$

for the chosen coordinates of the vertex $\mathbf{v_1}$. Now equation (14) for k=N-1 can be written as:

$$\mathbf{v_N} - \mathbf{v_1} + \mathbf{\Delta v_N} = \sum_{k=1..N} \mathbf{\Delta v_k}.$$

From equation (4) it is seen that both the x and y components of the right hand side of the above equation are zero valued. Thus flexibility of the juncture depends upon the value of the z component being 0; ie.

$$\sum_{k=1..N} \Delta z_k = 0. \qquad (16)$$

In the following paragraphs, equation (16) is shown to be true for each of the suspension types for specific choices of $\mathbf{v_1}$ and of signs in equations (15).

Turning to the suspension type **I-OEE** as defined by equations (3) we can now state the following theorem:

**Theorem I.** Prismatic junctures of Type I-OEE are flexible.

Proof: By direct substitution from equations (3) into equations (15) we have the following:

$$\Delta x_{k+M} = -\Delta x_k \text{ and}$$
$$\Delta y_{k+M} = \Delta y_k,$$

thus by taking

$$\Delta z_{k+M} = -\Delta z_k \text{ for } k = 1..M,$$

it is seen that equation (16) is true for any choice of the vertex $\mathbf{v_1}$. ∎

However, it is worth noting that by defining $\mathbf{v_1}$ as follows:

$$x_1 = \frac{-\sum_{k=1..M} \Delta x_k}{2},$$

$$y_1 = 0 \text{ and}$$

$$z_1 = \frac{-\sum_{k=1..M} \Delta z_k}{2},$$

we have axially symmetric coordinates for the juncture that are characterized by the following:

$$x_{k+M} = -x_k,$$
$$y_{k+M} = y_k \text{ and}$$
$$z_{k+M} = -z_k \text{ for } k = 1..M.$$

Thus it is seen that the axial symmetry that is found in Bricard octahedra of the first type is carried over into the I-OEE unbounded suspensions.

For suspension type **II-AEE** as defined by equations (8) we state the following theorem:

**Theorem II.** Prismatic junctures of Type II-AEE are flexible.

Proof: By direct substitution from equations (8) into equations (15) we have the following:

$$\Delta x_{N-k+1} = -\Delta x_k \text{ and}$$
$$\Delta y_{N-k+1} = \Delta y_k,$$

thus by taking

$$\Delta z_{N-k+1} = -\Delta z_k \text{ for } k = 1..M,$$

it is seen that equation (16) is true for any choice of the vertex $\mathbf{v_1}$. ∎

Further by defining $\mathbf{v_1}$ as follows:

$$x_1 = 0,$$
$$y_1 = 0 \text{ and}$$
$$z_1 = \frac{z_M}{2} - \sum_{k=1..M} \Delta z_k,$$

we have planar symmetric coordinates (as are found in Bricard octahedra of the second type) for the juncture that are characterized by the following:

$$x_{N-k+2} = x_k,$$
$$y_{N-k+2} = -y_k \text{ and}$$
$$z_{N-k+2} = -z_k \text{ for } k = 2..M.$$

For the suspension type **II-OEE** as defined by equations (11) we state the following theorem:

**Theorem III.** Prismatic junctures of Type II-OEE are flexible.

Proof: By direct substitution from equations (11) into equations (15) we have the following:
$$\Delta x_{k+M} = \Delta x_k \text{ and}$$
$$\Delta y_{k+M} = \Delta y_k,$$

thus by taking

$$\Delta z_{k+M} = -\Delta z_k \text{ for } k = 1..M,$$

it is seen that equation (16) is true for any choice of the vertex $\mathbf{v_1}$. ∎

Further by defining $\mathbf{v_1}$ as follows:

$$x_1 = 0,$$
$$y_1 = 0 \text{ and}$$
$$z_1 = \frac{-\sum_{k=1..M} \Delta z_k}{2},$$

we have another instance of planar symmetric coordinates for the juncture that are characterized by:

$$x_{k+M} = -x_k,$$
$$y_{k+M} = -y_k \text{ and}$$
$$z_{k+M} = z_k \text{ for } k = 2..M.$$

For the suspension type **III-OAE** as defined by equations (12) we state the following theorem:

**Theorem VI.** Prismatic junctures of Type III-OAE are flexible.

Proof: By direct substitution from equations (12) we have the following:

$$\Delta x_k = \frac{\pm L_k (\cos B - \cos \beta)}{2 \cos \theta} \text{ and}$$

$$\Delta y_k = \frac{\pm L_k (\cos B + \cos \beta)}{2 \cos \theta} \text{ for } k = 1..N.$$

Substituion of these terms into the equation for $\Delta z_k$ leads to:

$$\sum_{k=1..N} \Delta z_k = \left( \sum_{k=1..N} \pm L_k \right) \sqrt{1 - \frac{\cos^2 B + \cos^2 \beta}{2 \cos^2 \theta}}$$

The edge length relationships found in equations (13) provide the basis for the choice of signs for which the terms in the above summation vanish, thereby satisfying equation (16) for any choice of of vertex $\mathbf{v_1}$. ∎

We conclude with the following theorem related to the nature of the flexibility of the unbounded flexibile suspensions:

**Theorem VI.** All dihedral angles at the juncture and at the parallel edges of unbounded flexibile suspensions are non-constant under flexion.

Proof: It is easy to show that

$$\cos 2\theta = \cos \beta_k \cos B_k + \sin \beta_k \sin B_k \cos \varepsilon_k,$$

for all of the dihedral angles $\varepsilon_k$ for k=1..N at the juncture. Thus the dihedral angles $\varepsilon_k$ vary continuously as a function of the variable of flexion $\theta$. The dihedral angles $\delta_k$ along the parallel edges satisfy the "equation of the tetrahedral angle" [1; Eq. 1] and are shown there to be non-constant functions of the dihedral angles $\varepsilon_k$ at the juncture edges. ∎

## 4. Construction of Prismatic Flexible Polyhedra

This section describes how multiple unbounded suspensions are incorporated into flexible prismatic polyhedra that are comprised of annular segments, that are joined endwise at junctures, that are polygonal in cross section and each of which contains N (an even number ≥4) trapezoidal faces. These prismatic polyhedra are of genus 0 and 1 and are flexible by virtue of their ability to undergo changes in all of their dihedral angles while the geometric characteristics of their faces remain unchanged.

Genus 0 prismatic polyhedra have two or more segments that are joined endwise at junctures to form an open ended structure with unbounded segments at each end. The unbounded suspensions described in the previous section are examples of the simplest form of prismatic polyhedra; one having two segments. In general we count the two unbounded segment ends as two vertexes so that there are N(J-1)+2 vertexes, N(2J-1) edges and JN faces for J>1 segments; thus the genus 0 Euler equation (faces+vertexes=edges-2) holds. Genus 1 prismatic polyhedra have an even number of four or more segments that are joined endwise to form a torus; for J segments there are JN vertexes, 2JN edges and JN faces; the genus 1 Euler equation (faces+vertexes=edges) holds.

The flexibility of these prismatic polyhedra is achieved by constructions in which the geometry at each juncture of the respective segments satisfies the criteria found in equations (3), (8), (11) or (12). To facilitate the discussion we introduce some nomenclature and notation to describe polyhedra that are comprised of annular segments. Each segment is identified by an index, j say, that is associated with the face angles of the trapezoidal faces, the parallel edges of the faces and the edges at one of the junctures. The notation and indexing is illustrated in Fig. 5 for the k-th trapezoid of the j-th segment where k=1..N (the number of faces in a segment) and j=1..J (the number of segments in the polyhedron).

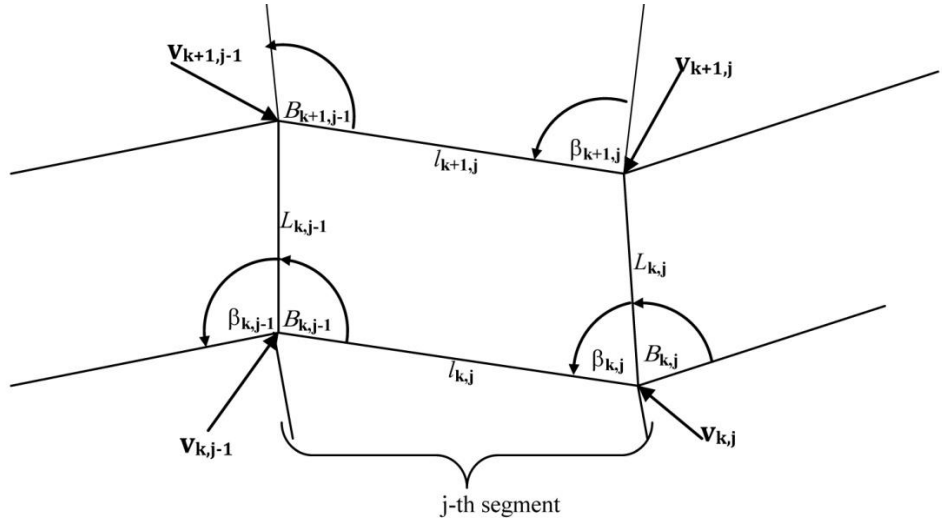

**Figure 5.** k-th Trapezoidal Face, j-th Annular Segment.

It is convenient to refer to segments and junctures by the notation $\mathsf{S}_j$ and $\mathsf{J}_j$ respectively. For both genus 0 and genus 1 there are J segments $\{\mathsf{S}_1,\mathsf{S}_2...\mathsf{S}_J\}$; for genus 0 polyhedra there are J-1 junctures $\{\mathsf{J}_1,\mathsf{J}_2...\mathsf{J}_{J-1}\}$; for genus 1 there are J junctures $\{\mathsf{J}_1,\mathsf{J}_2...\mathsf{J}_J\}$. Generally the juncture $\mathsf{J}_j$ is associated with the vertexes $\{\mathsf{v}_{k,j}; k=1..N\}$ while $\mathsf{S}_j$ is the segment that is terminated by the junctures $J_{j-1}$ and $J_j$. For genus 0 polyhedra the junctures $J_0$ and $J_J$ are effectively at infinity. For genus 1 the two junctures are the same, representing the same vertexes and parameters.

In this notation parameters at the juncture $J_j$, as shown in Fig. 2, are given by the set $\{\beta_{k,j}, B_{k,j}, L_{k,j}, k=1..N\}$. The notation for the first (j=1) and last (j=J) segments is exceptional; for genus 0 polyhedra the edge lengths $l_{k,1}$ and $l_{k,J}$ are of infinite extent and the vertexes $\mathbf{v}_{k,0}$ and $\mathbf{v}_{k,J}$ are at infinity. In a genus 1 polyhedra the edge lengths and adjacent angles are treated cyclically with the definition $\beta_{k,0}=\beta_{k,N}$ for example. For either case indexing with k is cyclic with $(N + 1) \equiv 1$.

Two "minimum lengths", $N_j$ and $M_j$, are associated with the juncture $J_j$ and are defined to be the smallest distances from the vertex $\mathbf{v_{1,j}}$ (taken along the trapezoidal edges) at which the cross section of the annular segments can be drawn and not intersect with the edges on the juncture. These lengths are defined in terms of the individual vertex offsets from vertex $\mathbf{v_{1,j}}$ ;

$$\nu_{k,j} = \nu_{k-1,j} - L_{k,j} \cos \beta_{k,j} \text{ and}$$
$$\mu_{k,j} = \mu_{k-1,j} + L_{k,j} \cos B_{k,j} \; k = 2..N \text{ with}$$
$$\nu_{1,j} = \mu_{1,j} = 0,$$

and by the expressions:

$$N_j = \min\{0, \nu_{k,j}, k = 1..(N-1)\} \text{ and}$$
$$M_j = \max\{0, \mu_{k,j}, k = 1..(N-1)\}.$$

The individual edge lengths $l_{k,j}$ of the segment $S_j$ are defined in terms of these minimum lengths and offset lengths and by a non-zero segment length $S_j$:

$$l_{k,j} = S_j - \mu_{k,j-1} - \nu_{k,j} + N_j + M_{j-1} \text{ for } k = 1..N. \tag{17}$$

The construction that is described in the following progresses by recursively adding segments to existing segments; at the j-th stage the segment $S_j$ is added to segment $S_{j-1}$ at juncture $J_{j-1}$. To initiate the recursion the left hand branch of the unbounded polyhedral suspension is treated as the first segment $S_1$, a full set of parameters that is consistent with the type of interest is assigned to juncture $J_1$ and the vertexes are positioned and edges are oriented as described in Sec. 3. The second segment $S_2$ and any subsequent segments are appended as described in the remainder of this section. Addition of the segment $S_j$ to segment $S_{j-1}$ at the juncture $J_{j-1}$ for j=2..J is accomplished by a specification of the segment length $S_j$ and of the definition of parameters at the juncture $J_j$. We limit our choice of parameters to the four cases represented by Fig. 1 and Fig. 3. It is easy to see that the orientation vectors associated with the alternate junctures (Figs. 3A, 3B, 3C) are defined by $(-\mathbf{w},\mathbf{u})$, $(-\mathbf{w},-\mathbf{u})$ and $(-\mathbf{u}, \mathbf{w})$ respectively.

For a particular choice of segment orientation it is possible to compute, from equation (17), individual edge lengths such that the vertexes of the new juncture are given by:

$$\mathbf{v_{k,j}} = \mathbf{v_{k,j-1}} + l_{k,j}\mathbf{w_j} \text{ for } k = 1..N. \tag{18}$$

Here $\mathbf{w_j}$ is a unit vector that defines the orientation of the parallel edges of segment $\mathbf{S_j}$. Choices of segment orientation are limited to one of the vectors $\mathbf{w}$, $-\mathbf{w}$, $\mathbf{u}$ and $-\mathbf{u}$, excluding the choice, and its negative, that was made at the previous junction. Note that for the first two segments $\mathbf{w_1} = -\mathbf{u}$ and $\mathbf{w_2} = \mathbf{w}$. By this approach flexible prismatic polyhedra of either genus 0 or 1 can be constructed. The latter when there are an even number of segments and the sum of segment lengths is zero when signs choosen for the orientation vectors are taken into account. However, motion of these polyhedra is essentially two dimensional as all orientation vectors lie in the same plane. A more general treatment in which the orientation of successive segments is not unconstrained is beyond the scope of this paper.

## 5. Conclusion

In this paper we have described four families of flexible prismatic polyhedra that can be developed from a consideration of Bricard octahedra and from flexible suspensions that result from generalizations of Bricard octahedra. These polyhedra are of genus 0 and genus 1, are of indefinite size, and consist of tubular segments formed from trapezoidal faces. We close with the following:

**Theorem VI.** Prismatic polyhedra, as described in the preceding section are flexible.

Proof: First, it is easy to see that the flexibility of junctures that was shown in Sec. 3 is applicable to each of alternate junctions as shown in Fig. 3. Second, geometric characteristics are preserved under the translation represented are by equation (18); thus the successive junctures are flexible. ∎

**Acknowledgements** The author wishes to thank his wife Verla Nelson for her encouragement and support.

**Gerald Nelson** is a retired software engineer; was employed by Honeywell, Inc. (when it was a Minnesota based company) and by MTS Systems Corp. of Eden Prairie, Minnesota. He has a Masters degree in Mathematics from the University of Minnesota.